\documentclass[12pt]{amsart}
\usepackage{amsfonts, amscd, amsmath}

\begin{document}

\newtheorem{theorem}{Theorem}[section]
\newtheorem{definition}[theorem]{Definition}
\newtheorem{proposition}[theorem]{Proposition}
\newtheorem{condition}[theorem]{Condition}
\newtheorem{lemma}[theorem]{Lemma}
\newtheorem{corollary}[theorem]{Corollary}
\theoremstyle{remark}
\newtheorem{example}[theorem]{Example}
\newtheorem{problem}[theorem]{Problem}

\def\BB{{\mathbb B}}
\def\CC{{\mathbb C}}
\def\EE{{\mathbb E}}
\def\FF{{\mathbb F}}
\def\RR{{\mathbb R}}
\def\ZZ{{\mathbb Z}}
\def\NN{{\mathbb N}}

\def\cB{{\mathcal B}}
\def\cD{{\mathcal D}}
\def\cF{{\mathcal F}}
\def\cH{{\mathcal H}}
\def\cJ{{\mathcal J}}
\def\cK{{\mathcal K}}
\def\cL{{\mathcal L}}
\def\cM{{\mathcal M}}
\def\cN{{\mathcal N}}
\def\cS{{\mathcal S}}
\def\cT{{\mathcal T}}

\def\Aut{\hbox{\rm Aut}\,}
\def\im{\hbox{\rm Im}\,}
\def\re{\hbox{\rm Re}\,}

\title[Characterization of the Hilbert ball]
{Characterization of the Hilbert ball \\
by its Automorphisms}
\author{Kang-Tae Kim and Daowei Ma}

\begin{abstract}
We show in this paper that every domain in a separable Hilbert
space, say $\cH$, which has a $C^2$ smooth strongly pseudoconvex
boundary point at which an automorphism orbit accumulates is
biholomorphic to the unit ball of $\cH$.  This is the complete
generalization of the Wong-Rosay theorem to a separable Hilbert
space of infinite dimension.  Our work here is an improvement from
the preceding work of Kim/Krantz [KIK] and subsequent
improvement of Byun/Gaussier/Kim [BGK] in the infinite dimensions.
\end{abstract}

\thanks{The first named author's research is supported in part by
The Grant KRF-2002-070-C00005 of The Korea Research Foundation.}

\subjclass[2000]{Primary: 32M05}
\address{Department of Mathematics, Pohang University of Science and
Technology, Pohang 790-784 Korea}
\email{kimkt@postech.edu}
\address{Department of Mathematics,
Wichita State University, Wichita, KS 67260-0033 U.S.A.}
\email{dma@math.twsu.edu}

\maketitle

\section{Introduction}

The primary goal of this article is to establish the following theorem, which gives a full generalization, to a separable Hilbert space of infinite dimension, of the Wong-Rosay Theorem of finite dimension.

\begin{theorem} \label{T:main}
If a domain $\Omega$ in a separable Hilbert space $\cH$ admits a $C^2$ strongly pseudoconvex boundary point at which a holomorphic automorphism orbit accumulates, then $\Omega$ is biholomorphic to the open unit ball in $\cH$.
\end{theorem}

Since there are many subtleties in setting up the necessary
terminology in the infinite dimensions, we shall present the
precise definitions in the next section.

There have been several important contributions by several authors
concerning this line of research. Chronologically speaking, Wong
[WON] proved in 1977 the above theorem in $\CC^n$ with the
assumption that the domain $\Omega$ is bounded and strongly
pseudoconvex at every boundary point.  Then, Rosay [ROS] improved
it in 1979, using holomorphic peak functions, that the theorem
holds if the domain is bounded and an automorphism orbit
accumulation point is strongly pseudoconvex.  Much later in 1995,
Efimov [EFI] removed the boundedness assumption from $\Omega$.
That argument is now well set up using Sibony's analysis of
plurisubharmonic peak functions. See [BER], [GAU] and [BGK], for
details. For the infinite dimension, Kim and Krantz [KIK] in 2000
proved the above theorem with an extra assumption that $\Omega$ is
bounded and convex. They needed convexity since they were relying
upon a weak-normal family argument which they developed. Then,
developing an infinite dimensional version of Sibony's analysis on
plurisubharmonic peak functions, Byun, Gaussier and Kim ([BGK] in
2002) removed the boundedness assumption from the theorem of Kim
and Krantz.  In this article, we remove the convexity assumption
from the theorem of Byun-Gaussier-Kim, thus arriving at the
optimal version of the theorem of this type. The crux of the proof
uses a new method, which concerns a principle of strong
convergence for certain holomorphic mappings of the infinite
dimensional Hilbert space. This new convergence argument seems
worth exploring further, with a separate interest. Finally, it is
worth noting that a manifold version of Wong-Rosay theorem have
been studied also. See [MAK] for instance.  Now the most general
version is known, and is due to Gaussier, Kim and Krantz ([GKK]).
We also present the Hilbert manifold version in this article.

The rest of the paper is organized as follows.  Since the proof
uses the ideas developed by Kim and Krantz [KIK] and
then the localization methods introduced in [BGK], we shall
introduce the outline of their methods shortly after the notation
and basic terminology are introduced.  Then we shall present our
methods leading to the strong convergence of the scaling sequence
in the separable Hilbert space and to the proof of the main
theorem.

\section{Terminology}

We introduce in this section the concepts of smoothness of mappings
of infinite dimensional spaces and the strong pseudoconvexity.
Further details in great generality can be found in the books of
Mujica [MUJ] and Dineen [DIN], for instance.

Let $\EE$ and $\FF$ be Banach spaces and let $\Omega$ be an open
subset of $\EE$.  Let $u:\Omega \to \FF$ be a $C^\infty$ smooth
mapping. Then for each point $p \in \Omega$ and vectors $v_1,
\ldots, v_k \in \EE$, we may define inductively the derivative
$d^k u$ of order $k$ as follows:
\begin{eqnarray*}
du (p; v_1) & = & \lim_{\RR \ni r \to 0} \frac1r (u(p+rv_1) - u(p)) \\
d^2 u (p; v_1, v_2) & = & \lim_{\RR \ni r \to 0} \frac1r
(du(p+rv_2; v_1) - du(p; v_1)) \\
& \vdots &
\end{eqnarray*}
\begin{multline*}
d^k u (p; v_1, \ldots, v_k) = \\
\lim_{\RR \ni r \to 0} \frac1r (d^{k-1} u (p+rv_k; v_1, \ldots,
v_{k-1}) - d^{k-1} u (p; v_1, \ldots, v_{k-1}))
\end{multline*}
Notice that these derivatives are symmetric multi-linear over
$\RR$.  If one so prefer, these formulae can be used to define
$C^k$ smoothness, requiring in that case that the corresponding
derivatives are continuous multi-linear tensors.

Then the complex differentials can be defined accordingly:
\begin{eqnarray*}
\partial u (p; v) & = & \frac12 (du(q;v) - i\, du(q; iv)) \\
\bar\partial u (p; v) & = & \frac12 (du(q;v) + i\, du(q; iv)).
\end{eqnarray*}

We are now able to introduce the concept of holomorphic maps and
the remaining terminology thereof.  First, by a {\it holomorphic
mapping} we mean a $C^1$ smooth map that is annihilated by the
$\bar\partial$ operator.  See [MUJ] for equivalent definitions.

A {\it domain} in this article is an open connected subset of a
Banach space.  An {\it automorphism} of a domain $\Omega$ is a
bijective holomorphic mapping of $\Omega$ with its inverse
holomorphic. The {\it automorphism group} $\Aut (\Omega)$ of a
domain $\Omega$ is the group of all automorphisms of $\Omega$. An
{\it automorphism orbit} is a set of the form $Aut(\Omega)q =
\{\varphi(q) \mid \varphi \in \Aut (\Omega)\}$, where $q\in
\Omega$. Thus, a boundary point $p$ of $\Omega$ is said to be an
{\it orbit accumulation point}, if $p$ is an accumulation point of
an automorphism orbit, i.e., if there is a sequence $\{\varphi_j
\}\subset \Aut (\Omega)$ and a point $q\in \Omega$ such that
$\displaystyle{\lim_{j\to \infty} \|\varphi_j (q)-p\| = 0}$.

We now introduce the concept of strong pseudoconvexity.  Let
$\Omega$ be a domain in a Banach space $\EE$.  We say that
$\Omega$ is strongly pseudoconvex at a boundary point $p \in
\partial\Omega$ if there is an open neighborhood $U$ of $p$ and a
$C^2$ smooth local defining function $\rho:U \to \RR$ satisfying
the following properties:
\begin{itemize}
\item[(\romannumeral1)] $\Omega \cap U = \{z \in U \mid \rho (z) <
0\}$.
\item[(\romannumeral2)] $\partial\Omega \cap U = \{z \in U \mid
\rho (z) = 0\}$.
\item[(\romannumeral3)] $d\rho (q; \cdot)$ is a non-zero
functional for every $q \in \partial\Omega \cap U$.
\item[(\romannumeral4)] There exists a constant $C>0$ such that
$\partial\bar\partial\rho (p; v, v) \ge C\|v\|^2$ for every $v \in
\EE$ satisfying $\partial\rho (p; v)=0$.
\end{itemize}

As in the finite dimensional case, we call $\partial\bar\partial
\rho$ the {\it Levi form} of $\rho$.

\section{Localization and Pluri-subharmonic Peak Functions}
\label{Localization}

Let $\Omega$ be a domain in a Banach space $\EE$, and let $p \in
\partial\Omega$.  A continuous function $\psi:\Omega \to \RR$ is
said to be {\it pluri-subharmonic} if it is subharmonic along
every complex affine line in $\Omega$.  A {\it
pluri-subharmonic peak function} at $p$ of $\Omega$ is a
pluri-subharmonic function $\psi_p : U \to \RR$ defined on an open
neighborhood $U$ of the closure $\overline{\Omega}$ of $\Omega$
satisfying the following two conditions:
\begin{itemize}
\item[(1)] $h(p)=0$, and $h(z)<0$ for every $z \in
\overline{\Omega} \setminus \{p\}$.
\item[(2)] The sets $V_m := \{z \in \overline{\Omega} \mid
h(z)>-1/m\}$, where $m=1,2,\ldots$, form a neighborhood basis at
$p$ in $\overline{\Omega}$.
\end{itemize}
Unlike the finite dimensional cases, the second condition is
essential for the definition in the infinite dimensions.  On the
other hand, notice that every strongly pseudoconvex boundary point
admits a pluri-subharmonic peak function for $\Omega$.

Following the work of Sibony [SIB], several investigations from
the articles of Efimov [EFI], Berteloot [BER], Gaussier [GAU] and
Byun- Gaussier-Kim [BGK] have been made.  We exploit some of them
which pertain to the localization and hyperbolicity.

\begin{theorem} {\rm (cf. p.\ 588, [BGK])} Let $\Omega$ be a domain in a
Banach space $\EE$ with a $C^2$ smooth strongly pseudoconvex
boundary point at which an automorphism orbit accumulates.  Then,
$\Omega$ is Kobayashi hyperbolic.
\end{theorem}

\begin{theorem} {\rm (cf. p.\ 588, [BGK])} Let $\Omega$ be a domin in a
Banach space $\EE$ with a $C^2$ smooth strongly pseudoconvex
boundary point $p$ which admits a sequence $\varphi_j \in \Aut
(\Omega)$ of automorphisms and a point $q \in \Omega$ such that
$\lim_{j\to\infty} \varphi_j (q) = p$.  Then, for every Kobayashi
distance ball $B_\Omega^K (x;r)$ of radius $r$ centered at $x \in
\Omega$ and for every open neighborhood $U$ of $p$ there exists
$N>0$ such that $\varphi_j (B_\Omega^K (x;r)) \subset U$ for every
$j > N$.
\end{theorem}

We choose not to include any details of the proofs, in order to
avoid an excessive repetition with the references cited above.
However, we consider it appropriate to point out that the
localization method using pluri-subharmonic peak functions
initiated by Sibony seems indeed more effective than the
traditional localization arguments relying upon holomorphic peak
functions and normal family arguments.

\section{Scaling Maps and Weak Normal Family}

The contents of this section are mostly from the article of Kim
and Krantz [KIK].  The scaling
method introduced here has its finite dimensional origin in the
work of Pinchuk [PIN], Frankel [FRA], Kim [KIM] and
others.  Some other details are in [BGK].

\subsection{Pinchuk's scaling sequence}

Here, we introduce Pinchuk's scaling sequence.  We begin with some
notation. We choose an orthonormal basis $e_1, e_2, \ldots$ for a
separable Hilbert space $\cH$.  Then for each $z \in \cH$, we write
$$
z = \sum_{m=1}^\infty z_m e_m,
$$
and
$$
z' = \sum_{m=2}^\infty z_m e_m.
$$

Let $\Omega$ be a domain in a separable Hilbert space $\cH$ with a
$C^2$ strongly pseudoconvex boundary point $p$.  Then, there exist
an open neighborhood $U$ of $p$ in $\cH$ and an injective
holomorphic mapping $G:U \to G(U) \subset \cH$ such that the
following hold:
\begin{itemize}
\item[(A)] $G(p)=0$.
\item[(B)] The domain $\Omega_U := G(\Omega\cap U)=\left\{z \in G(U)\mid \re z_1 > \psi (\im z_1, z') \right\} $ is strictly convex.
Moreover, the function $\psi$ is strongly convex and vanishes
precisely to the second order at the origin.
\end{itemize}

Now, as in the hypothesis of Theorem \ref{T:main}, we work with
the assumption that there exist $q \in \Omega$ and a sequence
$\varphi_j \in \Aut (\Omega)$ such that
$\displaystyle{\lim_{j\to\infty} \|\varphi_j (q) - p\| = 0}$.  We
may choose a subsequence if necessary so that we have $\varphi_j
(q) \in U$ for every $j=1,2,\ldots$.  Let $q_j =
G(\varphi_j (q))$.  Choose now for each $j$ the point $p_j \in
\partial \Omega_U$ such that $p_j - q_j = r_j e_1$ for some $r_j >
0$.  Note that $p_j' = q_j'$.  Let us denote by $p_{j1} = \langle
p_j, e_1\rangle$.  Then we consider a complex affine linear
isomorphism $H_j : \cH \to \cH : z \mapsto w$ defined by
\begin{eqnarray*}
w_1 &=& e^{\theta_j} (z_1 - p_{j1}) + T_j (z' - p') \\
w' & = & z' - p'
\end{eqnarray*}
where the bounded linear functional $T_j : (e_1)^\perp \to \CC$ is
chosen for each $j$ to satisfy the following properties:
\begin{itemize}
\item[(C1)] $H_j (\Omega_U)$ is supported by the real hyperplane
defined by $\re w_1 = 0$ at the origin of $\cH$.
\item[(C2)] $0 \in \partial T_j(\Omega_U)$.
\item[(C3)] $T_j (q_j) = e^{\theta_j} r_j$.
\end{itemize}
Note also that $e^{\theta_j}$ can be chosen so that it converges
to 1 as $j \to \infty$.  See [KIK] for an explicit choice for
these maps and the values for $\theta_j$.

Now we define the Pinchuk scaling sequence.  Define the linear map
$L_j : \cH \to \cH$ by
$$
L_j (w) = \frac{w_1}{r_j} e_1 + \frac{1}{\sqrt{r_j}} w'.
$$
Then the Pinchuk scaling sequence is defined by the composition
$$
\omega_j = L_j \circ H_j \circ G \circ \varphi_j.
$$
This map is not well-defined on $\Omega$.  However, the
localization theorems in Section \ref{Localization} implies that
there exists an increasing and exhausting sequence of Kobayashi
distance open balls $B_\Omega^K (q, R_k)$ ($k=1,2,\ldots; R_1 <
R_2 < \cdots$) so that $\omega_j$ is a well-defined map of
$B_\Omega^K (q, R_k)$ whenever $j \ge k$.

\subsection{Weak Normal Family Theorems}
In the finite dimensions, the Pinchuk scaling sequence $\omega_j$
defines a normal family whose subsequential limits are holomorphic
embeddings of $\Omega$ into the ambient Euclidean space.  However,
it is not the case in the infinite dimensions.  In this section,
we introduce the concept of weak convergence of holomorphic
mappings that produces holomorphic limits.  This is again from
[KIK] and [BGK].

\begin{theorem} {\rm (Theorem 4.4 of [BGK])}
Let $\EE$ be a separable Banach space, and let $\FF$ a reflexive
Banach space.  Let $\Omega_1$ and $\Omega_2$ be domains in $\EE$
and $\FF$, respectively.  Assume further that $\Omega_2$ is
bounded. Then, for every sequence $h_j:\Omega_1 \to \Omega_2$ of
holomorphic mappings, there exist a subsequence $h_{j_k}$ and a
holomorphic mapping $\widehat h : \Omega_1 \to \FF$ such that, for
each $x \in \Omega_1$, the sequence $h_{j_k} (x)$ converges weakly
to $\widehat h (x)$.
\end{theorem}

Unlike the finite dimensions, this theorem is not so effective.
The weak limit $\widehat h$ is holomorphic, but not in general
injective.  Also, it is not even guaranteed at this point that
$\widehat h (\Omega_1)$ is contained in $\Omega_2$.  (Although, we
do have that $\widehat h (\Omega_1)$ is contained in the closed
convex hull of $\Omega_2$ due to the reflexivity of $\FF$.)
Nonetheless, this is about the best one can obtain from the
general theory.

Before proceeding further, we point out that one can modify the
Pinchuk scaling sequence $\omega_j$ so that the range becomes
bounded.  In [KIK], a method is described in detail to modify
$\omega_j$ by composing with an explicit linear fractional
transformation $\Psi$ so that the map $\Psi \circ \omega_j$ has
its image in $(1+\epsilon_j) \BB$ for each $j$.  Moreover, the
sequence of positive numbers $\epsilon_j$ converges monotonically
to zero.

From here on, we shall denote the map $\Psi \circ \omega_j$ by
$\tau_j$ for $j=1,2,\ldots$.

\subsection{Calibration of Derivatives and
Kim-Krantz Scaling Sequence}
In [KIK], a new method of modifying the scaling sequence has been
introduced.  The goal is to have a strong convergence of the
derivatives $d\tau_j (q; \cdot)$, as $j \to \infty$.  In order to
do this, Kim and Krantz have used Hilbert isometries to calibrate
each differential $d\tau_j (q; \cdot)$ as follows: first they
consider the new basis $d\tau_j (q; e_m)$ ($m=1,2,\ldots$) for the
separable Hilbert space $\cH$.  Then they apply the Gram-Schmidt
process to these vectors such as
$$
f_{jm} := d\tau_j (q; e_m) - \sum_{k=1}^{m-1} \frac{\langle
d\tau_j(q; e_k), f_{jk} \rangle}{\langle f_{jk}, f_{jk} \rangle}
f_{jk}.
$$
It turns out that the vectors $f_{jm}$ have their norms uniformly
bounded away from zero by a constant independent of $j$ and $m$.
Then they define the Hilbert space isometries $S_j : \cH \to \cH$
arising from the condition $S_j (f_{jm}/\|f_{jm}\|) = e_m$ for
every $j,m=1,2, \ldots$.  Finally, the Kim/Krantz scaling sequence
$$
\sigma_j := S_j \circ \tau_j : B_\Omega^K (q; R_j) \to
(1+\epsilon_j) \BB
$$
is introduced.

If we use the notation $\Sigma_n = \CC e_1 \oplus \ldots \oplus
\CC e_n$, we can observe at this point immediately that $d\sigma_j
(q; \Sigma_n) = \Sigma_n$ for every positive integer $n$. Moreover
the sequences $\|d\sigma_j (q; \cdot)\|$ and $\|d\sigma_j (q;
\cdot)^{-1}\|$ are both uniformly bounded.  (See Section 7 of
[KIK] for details.)

In summary, one obtains the following:

\begin{proposition}
Let $\Omega$ be a domain in a separable Hilbert space $\cH$ with a
$C^2$ smooth, strongly pseudoconvex boundary point $p$ at which an
automorphism orbit accumulates.  Then, there exist a point $q \in
\Omega$, a decreasing sequence $\epsilon_j$ of positive numbers
tending to zero, an increasing sequence $R_j$ tending to infinity,
and a sequence of holomorphic mappings $\sigma_j : B_\Omega^K (q;
R_j) \to (1+\epsilon_j) \BB$ such that
\begin{itemize}
\item[(\romannumeral1)] $\sigma_j$ converges weakly to a
holomorphic mapping $\sigma$ at every point of $\Omega$;
\item[(\romannumeral2)] $\sigma_j(q)=0$ and $\sigma(q)=0$;
\item[(\romannumeral3)] $d\sigma_j(q)$ and $d\sigma(q)$ are {\it
calibrated} in the sense that they map the flag subspace
$\Sigma_n=\CC e_1\oplus \cdots\oplus \CC e_n$ into $\Sigma_n$ for
each positive integer $n$;

and

\item[(\romannumeral4)] $d\sigma_j(q)$ converges to $d\sigma(q)$
on every $\Sigma_n$.
\end{itemize}
\end{proposition}

Notice that the arguments up to this point are sufficient to prove
the main theorem of [KIK].  The main theorem of [BGK] is also in
the same line but uses more modifications for the convergence of
the sequence of $\sigma_j^{-1}$ since $\Omega$ may be convex but
still unbounded.  We would like to remark that our proof, as one
can see in the subsequent section, goes around such difficulties
establishing directly the two facts: (1) $\sigma$ is injective,
and (2) $\sigma (\Omega_1) = \Omega_2$.

\section{Strong Convergence of the Scaling Sequence}

\subsection{Techniques for Strong Convergence Arguments}
We now demonstrate a new method of strong normal families in the
infinite dimensional Hilbert space.
We begin with an estimate on the Kobayashi metric and distance. From
here on, $d_M$ and $k_M$ will denote the Kobayashi distance
and metric of the complex manifold $M$, respectively.  Let $u:
[0,1)\to [0,\infty)$ be defined by $u(t)=(1/2)\ln[(1+t)/(1-t)]$,
so that $u(t)=d_\Delta(0,t)$ and $u^{-1}(s)=\tanh s$.

\begin{lemma} \label{L:esti}
Let $\Omega\subset H$ be a Kobayashi hyperbolic domain. For $q,
x\in \Omega$, denote by $a=d_\Omega(x,q)$.  If $\Omega'$ is a
subdomain of $\Omega$ such that $\Omega'\supset \{y\in \Omega:
d_\Omega(y,q)<b\}$, where $b>a$, then $d_{\Omega'}(x,q)\le
a/\tanh(b-a)$, $k_{\Omega'}(x, v)\le k_{\Omega}(x, v)/ \tanh(b-a)$
for $v\in H$.
\end{lemma}

\begin{proof}
We first prove the second inequality. Let $s=\tanh(b-a)$ and
$\epsilon>0$. Then, there exists a holomorphic map $f: \Delta\to
\Omega$ such that $f(0)=x$ and $f'(0)=v/(k_\Omega(x,v)+\epsilon)$.
If $\zeta\in \Delta(0,s)$, then
\begin{eqnarray*}
d_\Omega(q,f(\zeta)) & \le & d_\Omega(q,x)+d_\Omega(x,f(\zeta))\\
& = & a+d_\Omega(f(0),f(\zeta)) \\
& \le & a+d_\Delta(0,\zeta) \\
& < & a+(b-a)=b.
\end{eqnarray*}
So $f(\Delta(0,s))\subset \Omega'$. Define $g:\Delta\to\Omega'$ by
$g(\zeta)=f(s\zeta)$. Then we have $g(0)=x$ and
$g'(0)=sf'(0)=sv/(k_\Omega(x,v)+\epsilon)$. Thus, it holds that
$k_{\Omega'}(x,v)\le (k_\Omega(x,v)+\epsilon)/s$.  Since
$\epsilon$ can be arbitrarily small, we obtain that
$k_{\Omega'}(x,v)\le k_\Omega(x,v)/s$.

We now prove the first inequality. Let $\epsilon\in (0, b-a)$.
There exists a $C^1$ curve $z: [0,1]\to\Omega$ such that $z(0)=q$,
$z(1)=x$, and $\int_0^1k_\Omega(z(t),z'(t))\,dt<a+\epsilon$. It
follows that $d_\Omega(q,z(t))<a+\epsilon<b$ and $z(t)\in \Omega'$
for each $t\in [0,1]$. By the inequality that we proved in the
preceding paragraph, $k_{\Omega'}(z(t),z'(t))\le
k_{\Omega}(z(t),z'(t))/\tanh(b-a-\epsilon)$. Therefore,
\begin{eqnarray*}
d_{\Omega'}(x,q) & \le & \int_0^1k_{\Omega'}(z(t),z'(t))\,dt \\
& \le & \frac1{\tanh(b-a-\epsilon)}
\int_0^1k_\Omega(z(t),z'(t))\,dt.
\end{eqnarray*}
So we see that $d_{\Omega'} (x,q)<
(a+\epsilon)/\tanh(b-a-\epsilon)$.  Allowing $\epsilon$ tend to
$0$, we obtain the desired conclusion.
\end{proof}

The next lemma follows immediately by a standard normal family
argument.

\begin{lemma} \label{L:disc}
For each positive number $\epsilon<1$ there exists a constant
$\delta>0$ such that for each holomorphic function
$f:\Delta\to\Delta$ with $f(0)=0$ and $f'(0)>1-\delta$ it holds
that
\begin{equation*}
|f(z)-z|< \epsilon, \;\;\;\;\;\; \mbox{whenever}\; |z|\le 1-\epsilon.
\end{equation*}
\end{lemma}

\bigskip

We now present a crucial lemma, which establishes the strong
uniform convergence on exhausting open subsets for the Kim/Krantz
scaling sequence introduced in Section 4.3. We begin with some
notation. For bounded linear operators $S$ and $T$ of the Hilbert
space $\cH$, we use the standard notation $S \le T$ (or $T\ge S$,
equivalently) which means that $\langle (T-S)x, x\rangle \ge 0$
for each $x\in \cH$. We also use two more standard notation: $I$
for the identity map of $\cH$, and the notation $\BB$ for the open
unit ball in $\cH$.

\begin{lemma} \label{L:ball}
Let $\{a_j\}$ be a sequence of positive numbers with $a_j\to 0$
and let $g_j: \BB \to \BB$ be a sequence of holomorphic mappings
such that $g_j(0)=0$ and $dg_j(0)\ge (1-a_j)I$. Then the sequence
$g_j$ converges to $I$ uniformly on each $r\BB$ with $0<r<1$.
\end{lemma}

\begin{proof}
Fix $0<r<1$ and $0<\epsilon<1/8$. Let $\zeta$ be a unit vector in
the Hilbert space $\cH$.  Define a function $f_j:\Delta\to \Delta$
by $f_j(z)=\langle g_j(z\zeta), \zeta\rangle$. Then
$f'_j(0)=\langle dg_j(0)\zeta, \zeta\rangle\ge 1-a_j$. By Lemma
\ref{L:disc}, there is a positive integer $k=k(r,\epsilon)$ such
that $|f_j(z)-z|<\epsilon$ whenever $j\ge k$ and $|z|\le r$. Let
$h_j(z) := g_j(z\zeta)- f_j(z)\zeta$. Then $\langle
h_j(z),\zeta\rangle =0$. By Schwarz's Lemma, we have $|z|^2\ge
\|g_j(z\zeta)\|^2$.  This implies that
$$
|z|^2 \ge |f_j(z)|^2+\|h_j(z)\|^2\ge
(|z|-\epsilon)^2+\|h_j(z)\|^2.
$$
Consequently, we obtain $\|h_j(z)\|^2<2\epsilon-\epsilon^2$ and
$$
|g_j(z\zeta)-z\zeta|^2=\|h_j(z)\|^2+\|f_j(z)\zeta-z\zeta\|^2\le
(2\epsilon-\epsilon^2)+\epsilon^2=2\epsilon
$$
for $j\ge k$ and $|z|\le r$. Therefore, we see immediately that the
sequence $g_j$ converges to $I$ uniformly on each $r\BB$.
\end{proof}

We introduce one more technical lemma before we present the proof
of our main theorem.

\begin{lemma} \label{L:final}
Let $\psi_j: \BB \to \BB$ be a sequence of holomorphic mappings
such that the sequence $\psi_j$ converges to $I$ uniformly on
$r\BB$ for every $0<r<1$. Then for each $0<r<1$ there exists a
$k\ge 1$ such that $\psi_j$ is injective on $r\BB$ and
$\psi_j(\BB)\supset r\BB$ whenever $j\ge k$.
\end{lemma}

\begin{proof}
By the Cauchy estimates, we can deduce that $d\psi_j$ converges to
$I$ uniformly on each $r\BB$. Thus, for fixed constants $r>0$ and
$\epsilon<1$, it holds that
$$
\|\psi_j(x)-\psi_j(y)-(x-y)\|<\epsilon\|x-y\|
$$
for any $x, y\in rB$ and any sufficiently large $j$. Observe that
the injectivity of $\psi_j$ on $r\BB$ follows from this inequality
immediately.

Now, let $0<r<1$. Choose $0<\epsilon<1/8$ so that
$(1+2\epsilon)r<1$. Let $j$ be sufficiently large so that
$\psi=\psi_j$ satisfies $\|d\psi-I\|<\epsilon$ on
$(1+2\epsilon)r\BB$. Fix $x\in r\BB$. Let $y_0=x$ and let
$y_k=x+y_{k-1}-\psi(y_{k-1})$ for $k=1,2,\ldots$. Then we see that
\begin{eqnarray*}
\|\psi(y_k)-x\| &=& \|\psi(y_k)-\psi(y_{k-1})-(y_k-y_{k-1})\| \\
&\le& \epsilon\|y_k-y_{k-1}\| \\
&=& \epsilon\|\psi(y_{k-1})-x\|.
\end{eqnarray*}
It follows that $\|y_{k+1}-y_k\|=\|\psi(y_k)-x\|\le
\epsilon^k\|y_1-y_0\|$. Thus $\|\psi(y_k)-x\|\to 0$ and $\{y_k\}$
is a Cauchy sequence. The completeness of $\cH$ implies that $y_k$
converges to a certain $y$.  Moreover, it is obvious now that
$\psi(y)=x$.  The remaining assertion follows immediately.
\end{proof}
\bigskip

\subsection{Proof of Theorem \ref{T:main}}
By the hypothesis, we are given a point $q\in \Omega$ such that
the $Aut(\Omega)$-orbit of $q \in \Omega$ accumulates at a
strongly pseudoconvex boundary point $p$.

For $j=2, 3, \dots$ let $b_j=u(1-1/j)=(1/2)\ln(2j-1)$. Recall that
the following have been proved in [BGK] (Also see Proposition 4.2
of this article):
\begin{itemize}
\item[(\romannumeral1)] $\Omega$ is hyperbolic.
\item[(\romannumeral2)] There exist subdomains $\Omega_j\subset
\Omega$, injective holomorphic mappings $\sigma_j:\Omega_j\to \cH$
($j=2, 3, \dots$), and a holomorphic mapping $\sigma: \Omega\to
\cH$, satisfying the following.

\begin{itemize}
\item[(a)] $\displaystyle{\bigcup_{j=1}^\infty \Omega_j=\Omega}$;
\smallskip
\item[(b)] $\Omega_j\supset \{x\in \Omega: d_\Omega(x,q)< b_j\}$;
\smallskip
\item[(c)] $\sigma_j(q)=0$, $\sigma(q)=0$;
\smallskip
\item[(d)] $d\sigma_j(q)$ and $d\sigma(q)$ are {\it calibrated} in
the sense that they map the flag subspace $\Sigma_n=\CC e_1\oplus
\cdots\oplus \CC e_n$ into $\Sigma_n$ for each positive integer
$n$;
\smallskip
\item[(e)] $(1-1/j)\BB\subset \sigma_j(\Omega_j)\subset(1+1/j)\BB$
and $\sigma(\Omega)\subset\overline \BB$;
\smallskip
\item[(f)] $\sigma_j$ converges weakly to $\sigma$ at every point
of $\Omega$;
\smallskip

and
\smallskip

\item[(g)] $d\sigma_j(q)$ converges to $d\sigma(q)$ on every
$\Sigma_n$.
\end{itemize}
\end{itemize}
\medskip
\noindent Our present goal is to show that $\sigma$ is a
biholomorphic mapping onto the open unit ball $\BB$ of the
separable Hilbert space $\cH$.

First we observe that $\sigma(\Omega)\subset \BB$ by the maximum
modulus principle. It follows by (b), (e) and Lemma \ref{L:esti}
that, for every $v\in \cH$,
\begin{equation}\label{eq:est}
 (1-1/j)k_\Omega(q, v)\le \|d\sigma_j(q)(v)\|
 \le (1+1/j)(1-1/j)^{-1}k_\Omega(q, v).
\end{equation}
With (g) and the fact that $d\sigma_j(q)$ are uniformly bounded,
we see that $\|d\sigma_j(q)(v)\|\to \|d\sigma(q)(v)\|$ for every
$v\in \cH$, as $j \to \infty$. It follows now that
\begin{equation}\label{eq:est0}
 \|d\sigma(q)(v)\|=k_\Omega(q, v).
\end{equation}

Consider $\sigma_j^{-1}: (1-1/j)\BB\to \Omega$. By (\ref{eq:est})
and (\ref{eq:est0}), it follows that
$$
\|d(\sigma\circ \sigma_j^{-1})(0)(v)\|\ge (1-1/j)(1+1/j)^{-1}\|v\|.
$$
Let $d(\sigma\circ \sigma_j^{-1})(0)=P_j U_j$ be the polar
decomposition of the invertible operator $d(\sigma\circ
\sigma_j^{-1})(0)$, where $P_j$ is positive and $U_j$ unitary.
Define a map $\tau_j: \BB\to\Omega$ by
$\tau_j(x)=\sigma_j^{-1}((1-1/j)U_j^{-1}x)$. Then $\sigma\circ
\tau_j: \BB\to \BB$, $\sigma\circ \tau_j(0)=0$.  Moreover, the
positive operator $d(\sigma\circ \tau_j)(0)=(1-1/j)P_j$ satisfies
$$
\|d(\sigma\circ \tau_j)(0)(v)\|\ge c_j\|v\|,
$$
where $c_j=(1-1/j)^2(1+1/j)^{-1}$.  It follows that $d(\sigma\circ
\tau_j)(0)\ge c_j I$, and that $c_j\to 1$. By Lemma~\ref{L:ball},
$\sigma\circ \tau_j$ converges to $I$ uniformly on $r\BB$ for
every $0<r<1$.

Fix $0<r<1$. By Lemma \ref{L:final},
$\sigma\circ\tau_j(\BB)\supset r\BB$ for sufficiently large $j$.
Hence $\sigma(\Omega)\supset r\BB$. Since this is true for each
$0<r<1$, we see that $\sigma(\Omega)=\BB$.

Fix $a>0$ and consider $Q_a=\{x\in \Omega: d_\Omega(x,q)<a\}$. Let
$r=(1+\tanh(a))/2$ and $t_j=\tanh(a/\tanh(b_j-a))$. If $j$ is
sufficiently large, then $b_j>a$ and $t_j(1+1/j)(1-1/j)^{-1}<r$.
By Lemma~\ref{L:esti}, $Q_a\subset \{x\in \Omega_j: d_{\Omega_j}(x,q)<a/\tanh(b_j-a)\}$. This, together with $\sigma_j(\Omega_j)\subset (1+1/j)\BB$, implies that $\sigma_j(Q_a)\subset
t_j(1+1/j)\BB$. It follows that
$$
\tau_j(r{\Bbb B})\supset\tau_j(t_j(1+1/j)(1-1/j)^{-1}\BB)
=\sigma_j^{-1}(t_j(1+1/j)\BB)\supset Q_a.
$$
For sufficiently large $j$, $\sigma\circ \tau_j$ is injective on
$r\BB$, hence $\sigma$ is injective on $\tau_j(r\BB)\supset Q_a$.
Since $\sigma$ is injective on $Q_a$ for every $a>0$, it must be
injective on $\Omega$. Therefore, $\sigma$ is an injective
holomorphic mapping from $\Omega$ onto $\BB$. It follows that
$\Omega$ is biholomorphic to $\BB$.  \qed

\section{Concluding Remarks}

Recall that the localization theorem is obtained from the
existence of pluri-subharmonic peak functions.  Since the
pluri-subharmonic functions are extremely flexible as far as the
extension properties are concerned, the whole argument of this
article is valid for the domains in a separable Hilbert manifold.
More precisely our main theorem extends to the following.

\begin{theorem}  Let $\Omega$ be a domain in a separable Hilbert manifold
$X$. If $\Omega$ admits an automorphism orbit
accumulating at a strongly pseudoconvex boundary point, then it is
biholomorphic to the unit ball in a separable Hilbert space $\cH$.
\end{theorem}

Notice that this is the infinite dimensional version of the main
theorem of [GKK].
\bigskip

The Wong-Rosay theorem in $\CC^n$ has been generalized to several
other domains that are not necessarily strongly pseudoconvex.
Better known theorems include [BEP], [KIM], [KIP] and [KKS].  They
characterized the Thullen domains and the polydiscs from Wong type
conditions on existence of boundary accumulating orbits.  However,
the authors do not know at this time of writing as how to
generalize these theorems to infinite dimensions.  Thus we close
this article posing the following question.

\begin{problem}
Formulate and prove the Wong-Rosay type characterization of the
unit ball in the space $c_0$ of complex sequences converging to
zero, or in the space $\ell^\infty$ of bounded complex sequences.
\end{problem}

\end{document}